\documentclass[11pt]{article}
\usepackage{eurosym}
\usepackage{amssymb}
\usepackage{amsmath}
\usepackage{amsfonts}
\usepackage{graphicx}
\usepackage{subfig}

\newtheorem{theorem}{Theorem}[section]

\newtheorem{corollary}[theorem]{Corollary}

\newtheorem{definition}[theorem]{Definition}

\begin{document}

\begin{center}
\vspace{3cm}

\textbf{\ On a Surface Pencil with a Common New Type of Special Surface
Curve in Galilean Space $G_{3}$}

\bigskip

Z\"{u}hal K\"{u}\c{c}\"{u}karslan Y\"{u}zba\c{s}\i , M\"{u}nevver Y\i ld\i
r\i m Y\i lmaz

F\i rat University, Faculty of Science, Department of Mathematics, 23119
Elazig / Turkey

zuhal2387@yahoo.com.tr,myildirim@firat.edu.tr

\bigskip
\end{center}

\textbf{Abstract: }In this study, we investigate a new type of a surface
curve called a new $D$-type special curve. Also, we show that this special
curve is more generally than a geodesic curve or an asymptotic curve. Then,%
\textbf{\ }we give the necessary and sufficient conditions for a curve to be
the new $D$-type special curve\ using Frenet frame in Galilean space. We
investigate some corollaries by taking account of a new $D$-type special
curve as a helix, a salkowski and an anti-salkowski. After all, for the sake
of visualizing of this study, we plot some examples for this surface pencil
(i.e. surface family).

\textbf{Key words: }New $D$-type special curve, Geodesic curve,
Isoparametric curve Parametric surface, Galilean space.

\section{Introduction}

\label{intro} One of the important problems arising when studying geometry
is to determine its properties by using physical,computational and
experimental methods. For this aim, many researchers have been focused on
the curve and surface theory because of having many applications to that of
various branch of science and engineering. As far as we know, the
Frenet-Frame based system is commonly used in physics as well confined
particle motion around the design orbit. The coordinate system around the
design orbit is called \ the Serret-Frenet frame and this \ frame can be
achieved using some special functions related to Hamiltonian, \cite{hwa,lee}.

On the other hand Galilean geometry is one of the real Cayley-Klein \
Geometries to that of motions are the Galilean transformations of classical
kinematics, \cite{Yaglom}. The decades have witnessed a rapid increase in
study of Galilean and Pseudo-Galilean space in \cite{Dede,dede2,alper,zuhal}.

Besides all studies mentioned above parametric of a surface pencil with a
common spatial geodesic, asymptotic and geometric applications for computer
science is have a great importance for those who study multidisciplinary
science and searching relations between theoretical and applied methodology,
in \cite{kasap,kasap2,li,wang}.

Roughly speaking, this work serves the purpose of defining a new type
surface curve which called a new D-type curve using Frenet-frame in Galilean
space. Furthermore, a $D-$type curve \ firstly defined as $\left\langle \eta
_{1},E_{0}\right\rangle =\lambda =const.$ in \cite{kaya}, then we introduce
a new D-type curve given by $\left\langle \eta _{1},E_{0}\mathbf{\wedge }%
t\right\rangle =\lambda =const.$ based on the first definition. We show that
this new $D-$type curve is more general than a geodesic or asymptotic curve.

This study is organized as follows: In preliminary part, we give Galilean
space $G_{3}$ and give some basic definitions and concepts of it. Then, we
define a new $D$-type special curve of this space. Following section is
devoted to surfaces with a common new type special curve in $G_{3}$ . We
give some characterizations for this new type curve. Notice that new $D$%
-type special curves are also satisfies being helix,salkowski and
anti-salkowski curves. At the end of the study, we give some examples for
this surface pencil.

\section{Preliminaries}

The Galilean space $\mathbf{G}_{3}$ is a one of the real Cayley-Klein space,
which has the projective metric of signature $(0,0,+,+)$.The absolute figure
of the Galilean space consists of an ordered triple$\{\omega ,f,I\}$ in
which $\omega $ is the ideal (absolute) plane, $\ f$ \ is the line (absolute
line) in $\omega $ and $I$ is the fixed elliptic involution of points of $f$%
. For more properties of Galilean space can be found in \cite{ali,roc}.

A plane is said to be Euclidean if it contains $f$, otherwise it is said to
be isotropic. In the given affine coordinates, isotropic vectors are of the
form $(0,x_{2},x_{3})$, whereas Euclidean planes are of the form $x=const.$
The induced geometry of a Euclidean plane is Euclidean and of an isotropic
plane isotropic (i.e. 2-dimensional Galilean or flag-geometry).

\begin{definition}
Let  \textbf{$a$}$=$ $\left( a_{1},a_{2},a_{3}\right) $ and \textbf{$b$}$=$ $%
\left( b_{1},b_{2},b_{3}\right) $ be vectors in $\mathbf{G}_{3}.$ A vector 
\textbf{$a$} is said to be isotropic if $a_{1}=0$, otherwise it is said to
be non-isotropic. Then the Galilean scalar product of these vectors is given
by
\end{definition}

\begin{equation*}
\left\langle a,b\right\rangle =\left\{ 
\begin{array}{cc}
a_{1}b_{1}, & \text{if }a_{1}\neq 0\text{ or }b_{1}\neq 0 \\ 
a_{2}b_{2}+a_{3}b_{3}, & \text{if \ }a_{1}=0\text{ and }b_{1}=0%
\end{array}%
\right. ,
\end{equation*}
\cite{Milin} .

\begin{definition}
Let \textbf{$a$}$=$ $\left( a_{1},a_{2},a_{3}\right) $ and \textbf{$b$}$=$ $%
\left( b_{1},b_{2},b_{3}\right) $ be vectors in $\mathbf{G}_{3},$ the cross
product of the vectors \textbf{$a$ }and \textbf{$b$} is defined by 
\begin{equation*}
\ a\wedge b=\left\vert 
\begin{array}{ccc}
0 & e_{2} & e_{3} \\ 
a_{1} & a_{2} & a_{3} \\ 
b_{1} & b_{2} & b_{3}%
\end{array}%
\right\vert =(0,a_{3}b_{1}-a_{1}b_{3},a_{1}b_{2}-a_{2}b_{1}),
\end{equation*}
\cite{Milin}.
\end{definition}

\begin{definition}
If an admissible curve\ $r$ of the class $\mathbf{C}^{r}$ $(r>3)$ in $%
\mathbf{G}_{3},$ and parametrized by the invariant parameter $s$ , is given
by%
\begin{equation*}
r\left( s\right) =\left( s,f\left( s\right) ,g\left( s\right) \right) .
\end{equation*}%
then $s$ is a Galilean invariant of the arc length on $r$.
\end{definition}

The moving trihedron is written by 
\begin{eqnarray*}
t\left( s\right)  &=&r^{\prime }\left( s\right) =\left( 1,f^{\prime }\left(
s\right) ,g^{\prime }\left( s\right) \right) , \\
n\left( s\right)  &=&\frac{r^{\prime \prime }\left( s\right) }{\kappa \left(
s\right) }=\frac{1}{\kappa \left( s\right) }\left( 0,f^{\prime \prime
}\left( s\right) ,g^{\prime \prime }\left( s\right) \right) , \\
b\left( s\right)  &=&\frac{1}{\kappa \left( s\right) }\left( 0,-g^{\prime
\prime }\left( s\right) ,f^{\prime \prime }\left( s\right) \right) ,
\end{eqnarray*}%
where $t,n$ and $b$ are called the vectors of the tangent, principal normal
and binormal of $r\left( s\right) ,$ respectively, and the curvature $\kappa
\left( s\right) $ and the torsion $\tau \left( s\right) $ of the curve $r$
can be given by, respectively,%
\begin{eqnarray*}
\kappa \left( s\right)  &=&\sqrt{f^{\prime \prime }\left( s\right)
^{2}+g^{\prime \prime }\left( s\right) ^{2}}, \\
\tau \left( s\right)  &=&\frac{\det \left( r^{\prime }\left( s\right)
,r^{\prime \prime }\left( s\right) ,r^{\prime \prime \prime }\left( s\right)
\right) }{\kappa ^{2}\left( s\right) }.
\end{eqnarray*}%
Frenet formulas can be given as%
\begin{eqnarray}
t^{\prime } &=&\kappa n,  \label{1} \\
n^{\prime } &=&\tau b,  \notag \\
b^{\prime } &=&-\tau n,  \notag
\end{eqnarray}
\cite{pav}.

\begin{definition}
Let $\left\{ t,n,b\right\} $ be the Frenet frame of the differentiable curve
in $\mathbf{G}_{3}$. The equations 
%TCIMACRO{\TeXButton{1}{\eqref{1}} }%
%BeginExpansion
\eqref{1}
%EndExpansion
form a rotation motion with Darboux vector $w=\tau t+\kappa b.$ Also
momentum rotation vector holds the following conditions 
\begin{equation*}
t^{\prime }=w\mathbf{\wedge }t\mathbf{,}\text{ }n^{\prime }=w\mathbf{\wedge }%
n\mathbf{,}\text{ }b^{\prime }=w\mathbf{\wedge }b\mathbf{,}
\end{equation*}%
\cite{roc}
\end{definition}

Besides, the unit vector $E_{0}$ from Darboux vector is obtained by 
\begin{equation}
E_{0}=\frac{\tau t+\kappa b}{\left\vert \tau \right\vert },  \label{2}
\end{equation}%
\cite{kara}.

In the rest of paper, we shall suppose $\kappa\neq{0}$ and $\tau\neq{0}$ at
everywhere.

\begin{definition}
Let $\varphi $ is a surface in $\mathbf{G}_{3},$ the equation of a surface
in $\mathbf{G}_{3}$ can be expressed as the parametrization%
\begin{equation*}
\varphi \left( s,\upsilon \right) =\left( \varphi _{1}\left( s,\upsilon
\right) ,\varphi _{2}\left( s,\upsilon \right) ,\varphi _{3}\left(
s,\upsilon \right) \right) ,\text{ }s,\upsilon \in R,
\end{equation*}%
where $\varphi _{1}\left( s,\upsilon \right) ,\varphi _{2}\left( s,\upsilon
\right) $ and $\varphi _{3}\left( s,\upsilon \right) \in \mathbf{C}^{3},$ in 
\cite{roc}.
\end{definition}

Also, the isotropic normal vector is defined by 
\begin{equation*}
\eta \left( s,v\right) =\varphi _{s}\mathbf{\wedge } \varphi _{v},
\end{equation*}%
where $\varphi _{s}=\frac{\partial \varphi \left( s,v\right) }{\partial s}$
and $\varphi _{v}=\frac{\partial \varphi \left( s,v\right) }{\partial v}$.

\begin{definition}
Let $\varphi $ be a regular surface in $\mathbf{G}_{3}$ with the isotropic
surface normal $\eta \left( s,v\right) $ and $r(s)$ be an arc-length
parametrized curve on '$\varphi .$ If the following condition 
\begin{equation*}
\left\langle \eta _{1},E_{0}\mathbf{\wedge }t\right\rangle =\lambda =const.,
\end{equation*}%
is satisfied, then the curve $r(s)$ is said to be a new $D$-type special
curve on $\varphi ,$ where $E_{0},$ $t$ and $\eta _{1}$ is the unit darboux
vector$,$ tangent vector and unit surface normal along the curve $r$,
respectively.
\end{definition}

Considering Definition 2.6, if $\lambda =0$, then the surface normal $\eta
_{1}$ is orthogonal to the principal normal $n,$ i.e, the curve $r(s)$ is an
asymptotic curve on $\varphi .$ Similarly, if $\lambda =1$ then the surface
normal $\eta _{1}$ and the principal normal $\eta _{1}$ are linearly
dependent$,$ it means that the curve $r(s)$ is a geodesic curve on $\varphi
. $ Our studies show that the new $D$-type special curves contain both
geodesic and asymptotic curves, that is, the new $D$-type special curves are
more general then both curves.

\section{Surfaces with Common New Type Special Curve in Galilean Space $G_3$}

\label{sec:2} Let $\varphi =\varphi \left( s,v\right) $ be a parametric
surface on the arc-length parametrized curve $r(s)$ in $\mathbf{G}_{3}.$ The
surface is defined by%
\begin{eqnarray}
\varphi \left( s,v\right) &=&r(s)+[\alpha \left( s,v\right) t\left( s\right)
+\beta \left( s,v\right) n\left( s\right)  \label{3} \\
&&+\gamma \left( s,v\right) b\left( s\right) ],  \notag \\
L_{1} &\leq &s\leq L_{2}\text{ and }T_{1}\leq v\leq T_{2},  \notag
\end{eqnarray}%
where $\alpha \left( s,v\right) ,\beta \left( s,v\right) $ and $\gamma
\left( s,v\right) $ are smooth functions. $\alpha \left( s,v\right) ,\beta
\left( s,v\right) $ and $\gamma \left( s,v\right) $ are smooth functions and
their values indicate, respectively, the extension-like, flexion-like, and
retortion-like effects, by the point unit through the time $v$, starting
from $r(s),$ (see \cite{zuhal}).

Our starting point is to provide the necessary and sufficient conditions for
the given curve $r(s)$ to be a new $D-$type special curve on the surface $%
\varphi =\varphi \left( s,v\right) .$

The unit surface normal $\eta _{1}$ can be given%
\begin{equation*}
\eta _{1}=\cos \theta n+\sin \theta b,
\end{equation*}%
where $n$ and $b$ are the principal normal and binormal of $r\left( s\right)
,$ respectively. Now we give\ the necessary and sufficient conditions for an
isoparametric curve $r(s)$ to be a common special new $D-$type curve on $%
\varphi \left( s,v\right) .$

The isotropic normal $\eta \left( s,v\right) $ of the surface is given by 
\begin{equation}
\eta \left( s,v\right) =\varphi _{s}\mathbf{\wedge }\varphi _{v},  \label{4}
\end{equation}%
from 
%TCIMACRO{\TeXButton{3}{\eqref{3}} }%
%BeginExpansion
\eqref{3}
%EndExpansion
\begin{eqnarray*}
\varphi _{s} &=&\left( 1+\alpha _{s}\right) t+\left( k\alpha +\beta
_{s}-\tau \gamma \right) n+\left( \tau \beta +\gamma _{s}\right) b, \\
\varphi _{v} &=&\alpha _{v}t+\beta _{v}n+\gamma _{v}b.
\end{eqnarray*}%
Taking account 
%TCIMACRO{\TeXButton{4}{\eqref{4}}}%
%BeginExpansion
\eqref{4}%
%EndExpansion
$,$ $\eta \left( s,v\right) $ can be given as%
\begin{eqnarray*}
\eta \left( s,v\right)  &=&\left[ -\left( 1+\alpha _{s}\right) \gamma
_{v}+\left( \tau \beta +\gamma _{s}\right) \alpha _{v}\right] n \\
&&+\left[ \left( 1+\alpha _{s}\right) \beta _{v}-\left( k\alpha +\beta
_{s}-\tau \gamma \right) \alpha _{v}\right] b.
\end{eqnarray*}%
$.$

Let $r(s)$ be a curve on a surface $\varphi \left( s,v\right) $ in $\mathbf{G%
}_{3}.$ If $r(s)$ is isoparametric curve on this surface, then there exists
a parameter $v=v_{0}$ such that $r(s)=\varphi \left( s,v_{0}\right) ,$ that
is 
\begin{equation}
\alpha \left( s,v_{0}\right) =\beta \left( s,v_{0}\right) =\gamma \left(
s,v_{0}\right) =0.  \label{5}
\end{equation}%
From 
%TCIMACRO{\TeXButton{5}{\eqref{5}}}%
%BeginExpansion
\eqref{5}%
%EndExpansion
$,$ we get 
\begin{eqnarray}
\eta \left( s,v_{0}\right) &=&\left[ -\left( 1+\alpha _{s}\right) \gamma
_{v}+\gamma _{s}\alpha _{v}\right] n  \label{6} \\
&&+\left[ \left( 1+\alpha _{s}\right) \beta _{v}-\beta _{s}\alpha _{v}\right]
b.  \notag
\end{eqnarray}
From 
%TCIMACRO{\TeXButton{6}{\eqref{6}}}%
%BeginExpansion
\eqref{6}%
%EndExpansion
, we can write 
\begin{eqnarray*}
\varphi _{1}\left( s,v_{0}\right) &=&0, \\
\varphi _{2}\left( s,v_{0}\right) &=&\left[ -\left( 1+\alpha _{s}\right)
\gamma _{v}+\left( \tau \beta +\gamma _{s}\right) \alpha _{v}\right] , \\
\varphi _{3}\left( s,v_{0}\right) &=&\left[ \left( 1+\alpha _{s}\right)
\beta _{v}-\left( k\alpha +\beta _{s}-\tau \gamma \right) \alpha _{v}\right]
.
\end{eqnarray*}

Since $\eta _{1}$ is parallel to $\eta \left( s,v_{0}\right) ,$ then there
exists a function $\sigma (s)$ such that 
\begin{equation}
\varphi _{1}=0,\varphi _{2}=\sigma (s)\cos \theta ,\text{ }\varphi
_{3}=\sigma (s)\sin \theta .  \label{7}
\end{equation}%
Hence, the necessary and sufficient conditions for the surface $\varphi $ to
have the curve $r(s)$ as the new $D-$type special curves can be given with
the following theorem.

\begin{theorem}
Let $\varphi \left( s,v\right) $ be a surface having a curve $r(s)$ in $%
\mathbf{G}_{3}$ 
%TCIMACRO{\TeXButton{3}{\eqref{3}}}%
%BeginExpansion
\eqref{3}%
%EndExpansion
. The curve $r(s)$ is a new $D$-type special curve on a surface $\phi $ if
and only if \textbf{\ }%
\begin{equation*}
\alpha \left( s,v_{0}\right) =\beta \left( s,v_{0}\right) =\gamma \left(
s,v_{0}\right) =0,
\end{equation*}%
\begin{eqnarray}
\varphi _{1} &=&0,\text{ }\varphi _{2}=\lambda \frac{\left\vert \tau
\right\vert }{\kappa },  \label{8} \\
\varphi _{3} &=&\pm \sigma \sqrt{1-\left( \frac{\lambda }{\sigma }\frac{\tau 
}{\kappa }\right) ^{2}},  \notag
\end{eqnarray}%
satisfy, where $0\leq s\leq L$ and $0\leq v_{0}\leq T,\sigma (s)\neq 0,$ $%
\lambda $ is a real constant $\kappa $, $\tau $ are the curvature and the
torsion function of \ $r(s),$ respectively.
\end{theorem}

\textbf{Proof. }Let $r(s)$ be a special $D-$type curve on surface pencil $%
\varphi \left( s,v\right) .$ From Definition 2.6, we have 
\begin{equation}
\left\langle \eta _{1},E_{0}\mathbf{\wedge }t\right\rangle =\lambda ,
\label{9}
\end{equation}%
where $\lambda $ is a real constant. From 
%TCIMACRO{\TeXButton{9}{\eqref{9}}}%
%BeginExpansion
\eqref{9}%
%EndExpansion
, we get%
\begin{equation*}
\left\langle \eta _{1},n\right\rangle =\lambda \frac{\left\vert \tau
\right\vert }{\kappa }.
\end{equation*}%
By the taking account of $\eta _{\text{ }}$and $\eta _{1}$ are parallel to
each other, we can give

\begin{equation*}
\left\langle \varphi _{2}\left( s,v_{0}\right) n(s)+\varphi _{3}\left(
s,v_{0}\right) b(s),n(s)\right\rangle =\lambda \frac{\left\vert \tau
\right\vert }{\kappa },
\end{equation*}%
then, we obtain 
\begin{equation*}
\varphi _{2}\left( s,v_{0}\right) =\sigma (s)\cos \theta =\lambda \frac{%
\left\vert \tau \right\vert }{\kappa },
\end{equation*}%
Using 
%TCIMACRO{\TeXButton{7}{\eqref{7}}}%
%BeginExpansion
\eqref{7}%
%EndExpansion
, we get%
\begin{equation*}
\varphi _{3}\left( s,v_{0}\right) =\pm \sigma \sqrt{1-\left( \frac{\lambda }{%
\sigma }\frac{\tau }{\kappa }\right) ^{2}}.
\end{equation*}%
From 
%TCIMACRO{\TeXButton{3}{\eqref{3}}}%
%BeginExpansion
\eqref{3}%
%EndExpansion
, we have $\alpha \left( s,v_{0}\right) =\beta \left( s,v_{0}\right) =\gamma
\left( s,v_{0}\right) =0.$

Conversely, suppose that

\begin{equation*}
\alpha \left( s,v_{0}\right) =\beta \left( s,v_{0}\right) =\gamma \left(
s,v_{0}\right) =0,
\end{equation*}%
\begin{eqnarray*}
\varphi _{1} &=&0,\text{ }\varphi _{2}=\lambda \frac{\left\vert \tau
\right\vert }{\kappa }, \\
\varphi _{3} &=&\pm \sigma \sqrt{1-\left( \frac{\lambda }{\sigma }\frac{\tau 
}{\kappa }\right) ^{2}},
\end{eqnarray*}%
satisfy. Then 
%TCIMACRO{\TeXButton{3}{\eqref{3}} }%
%BeginExpansion
\eqref{3}
%EndExpansion
holds and for the surface normal along curve, we have

\begin{equation*}
\eta \left( s,v_{0}\right) =\lambda \frac{\left\vert \tau \right\vert }{%
\kappa }n+\pm \sigma \sqrt{1-\left( \frac{\lambda }{\sigma }\frac{\tau }{%
\kappa }\right) ^{2}}b
\end{equation*}%
and we get the following equations%
\begin{equation*}
\left\langle \eta \left( s,v_{0}\right) ,E_{0}\mathbf{\wedge }t\right\rangle
=\lambda =const.
\end{equation*}%
Because the vectors $\eta _{1}(s)$ and $\eta \left( s,v_{0}\right) $ are
parallel, we obtain 
\begin{equation*}
\left\langle \eta _{1},E_{0}\mathbf{\wedge }t\right\rangle =const.
\end{equation*}%
Hence the curve $r(s)$ is the new $D$-type special curve on surface pencil $%
\varphi \left( s,v\right) .$

From the above theorem, we have the following corollaries:

\begin{corollary}
\textbf{\ }Let the curve $r(s)$ be a new $D$-type special curve on the
surface $\varphi \left( s,v\right) .$ Then $r(s)$ is an isogeodesic curve on 
$\varphi \left( s,v\right) $ iff the following conditions are satisfied:%
\begin{equation*}
\alpha \left( s,v_{0}\right) =\beta \left( s,v_{0}\right) =\gamma \left(
s,v_{0}\right) =0,
\end{equation*}%
\begin{eqnarray}
\varphi _{1} &=&0,\text{ }\varphi _{2}=\frac{\left\vert \tau \right\vert }{%
\kappa },  \label{10} \\
\varphi _{3} &=&\pm \sigma \sqrt{1-\left( \frac{\tau }{\sigma \kappa }%
\right) ^{2}},  \notag
\end{eqnarray}%
where $0\leq s\leq L$ and $0\leq v_{0}\leq T,\sigma (s)\neq 0.$
\end{corollary}

\begin{corollary}
\textbf{\ }Let the curve $r(s)$ be a new $D$-type special curve on the
surface $\varphi \left( s,v\right) .$ Then $r(s)$ is an isoasymptotic curve
on $\varphi \left( s,v\right) $ iff the following conditions are satisfied:%
\begin{equation*}
\alpha \left( s,v_{0}\right) =\beta \left( s,v_{0}\right) =\gamma \left(
s,v_{0}\right) =0,
\end{equation*}%
\begin{eqnarray}
\varphi _{1} &=&0,\text{ }\varphi _{2}=0,  \label{11} \\
\varphi _{3} &=&\pm \sigma ,  \notag
\end{eqnarray}%
where $0\leq s\leq L$ and $0\leq v_{0}\leq T,\sigma (s)\neq 0.$
\end{corollary}

\begin{corollary}
Let the curve $r(s)$ be a new $D$-type special curve on $\varphi \left(
s,v\right) .$ Then $r(s)$ is a general helix on $\varphi \left( s,v\right) $
iff the following conditions are satisfied:%
\begin{equation*}
\alpha \left( s,v_{0}\right) =\beta \left( s,v_{0}\right) =\gamma \left(
s,v_{0}\right) =0,
\end{equation*}%
\begin{eqnarray}
\varphi _{1} &=&0,\text{ }\varphi _{2}=\lambda \mu ,  \label{12} \\
\varphi _{3} &=&\pm \sigma \sqrt{1-\left( \frac{\lambda }{\sigma }\mu
\right) ^{2}},  \notag
\end{eqnarray}%
where $0\leq s\leq L$ and $0\leq v_{0}\leq T,\sigma (s)\neq 0,$ $\lambda $
and $\frac{\left\vert \tau \right\vert }{\kappa }=\mu $ are real constants$.$
\end{corollary}

%\textbf{Corollary 3.4. }Let the curve $r(s)$ be a new $D$-type special curve
%on the surface $\varphi \left( s,v\right) .$ Then $r(s)$ is a planar curve
%on $\varphi \left( s,v\right) $ iff the following conditions are satisfied:%
%\begin{equation*}
%\alpha \left( s,v_{0}\right) =\beta \left( s,v_{0}\right) =\gamma \left(
%s,v_{0}\right) =0,
%\end{equation*}%
%\begin{eqnarray}
%\varphi _{1} &=&0,\text{ }\varphi _{2}=0,  \label{13} \\
%\varphi _{3} &=&\pm \sigma ,  \notag
%\end{eqnarray}%
%where $0\leq s\leq L$ and $0\leq v_{0}\leq T,\sigma (s)\neq 0.$

\begin{corollary}
\textbf{\ }Let the curve $r(s)$ be a new $D$-type special curve on the
surface $\varphi \left( s,v\right) .$ Then $r(s)$ is a Salkowski curve (or a
slant helices) on $\varphi \left( s,v\right) $ iff the following conditions
are satisfied:%
\begin{equation*}
\alpha \left( s,v_{0}\right) =\beta \left( s,v_{0}\right) =\gamma \left(
s,v_{0}\right) =0,
\end{equation*}%
\begin{eqnarray}
\varphi _{1} &=&0,\text{ }\varphi _{2}=\lambda \frac{\left\vert \tau
\right\vert }{\nu },  \label{14} \\
\varphi _{3} &=&\pm \sigma \sqrt{1-\left( \frac{\lambda }{\sigma }\frac{\tau 
}{\nu }\right) ^{2}},  \notag
\end{eqnarray}%
where $0\leq s\leq L$ and $0\leq v_{0}\leq T,\sigma (s)\neq 0,$ $\lambda $
and $\kappa =\nu $ are real constants and $\tau $ is non-constant.
\end{corollary}

\begin{corollary}
\textbf{\ }Let the curve $r(s)$ be a new $D$-type special curve on the
surface $\varphi \left( s,v\right) .$ Then $r(s)$ is an anti-Salkowski curve
(or a slant helices) on $\varphi \left( s,v\right) $ iff the following
conditions are satisfied:%
\begin{equation*}
\alpha \left( s,v_{0}\right) =\beta \left( s,v_{0}\right) =\gamma \left(
s,v_{0}\right) =0,
\end{equation*}%
\begin{eqnarray}
\varphi _{1} &=&0,\text{ }\varphi _{2}=\lambda \frac{\left\vert \xi
\right\vert }{\kappa },  \label{15} \\
\varphi _{3} &=&\pm \sigma \sqrt{1-\left( \frac{\lambda }{\sigma }\frac{\xi 
}{\kappa }\right) ^{2}},  \notag
\end{eqnarray}%
where $0\leq s\leq L$ and $0\leq v_{0}\leq T,\sigma (s)\neq 0,$ $\lambda $
and $\tau =\xi $ are real constants and $\kappa $ is non-constant $.$
\end{corollary}

Now, we can express the marching-scale functions $\alpha \left( s,v\right)
,\beta \left( s,v\right) $ and $\gamma \left( s,v\right) $ as the product of
two valued $C^{1}$ functions. Then we can give%
\begin{eqnarray}
\alpha \left( s,v\right) &=&l\left( s\right) X\left( v\right) ,  \notag \\
\beta \left( s,v\right) &=&m\left( s\right) Y\left( v\right) ,  \label{16} \\
\gamma \left( s,v\right) &=&n\left( s\right) Z\left( v\right) ,  \notag
\end{eqnarray}%
where $l$ $\left( s\right) ,m\left( s\right) ,n\left( s\right) $, $X\left(
v\right) ,Y\left( v\right) ,Z\left( v\right) $ are $C^{1}$ functions and $%
l\left( s\right) ,m\left( s\right) $ and $n\left( s\right) $ are not
identically zero.

Therefore, we can express the following corollary:

\begin{corollary}
\textbf{\ }The curve $r(s)$ is a new $D$-type special curve on the surface
pencil $\varphi \left( s,v\right) $ iff the following conditions are
satisfied:%
\begin{eqnarray}
X\left( v_{0}\right)  &=&Y\left( v_{0}\right) =Z\left( v_{0}\right) =0, 
\notag \\
m\left( s\right) \beta ^{\prime }\left( v_{0}\right)  &=&\pm \sigma \sqrt{%
1-\left( \frac{\lambda }{\sigma }\frac{\tau }{\kappa }\right) ^{2}},
\label{17} \\
-n\left( s\right) \gamma ^{\prime }\left( v_{0}\right)  &=&\lambda \frac{%
\left\vert \tau \right\vert }{\kappa },  \notag
\end{eqnarray}%
where $0\leq s\leq L$ and $0\leq v_{0}\leq T,\sigma (s)\neq 0,$ $\lambda $
is real constant and $\kappa $ and $\tau $ are the curvature and the torsion
functions of the curve $r(s)$ , respectively.
\end{corollary}

\bigskip

\textbf{Example 3.8}Let $r(s)$ be a general helix given by parametrization
in $\mathbf{G}_{3}$

\begin{equation*}
r\left( s\right) =\left( s,8\sqrt{\pi }FresnelS(\frac{1}{\sqrt{2\pi }}s),-8%
\sqrt{\pi }FresnelC(\frac{1}{\sqrt{2\pi }}s)\right) ,
\end{equation*}%
where $FresnelS(\gamma )=\int \sin (\frac{\pi \gamma ^{2}}{2})d\gamma $ and $%
FresnelC(\gamma )=\int \cos (\frac{\pi \gamma ^{2}}{2})d\gamma $, \cite{ali}%
. The plot of the curve is given by Fig 1a. It is easy to calculate that $%
\kappa =s$ and $\tau =\frac{1}{4}\kappa $, 
\begin{eqnarray*}
t &=&\left( 1,4\sin \frac{s^{2}}{8},-4\cos \frac{s^{2}}{8}\right) , \\
n &=&\left( 0,\cos \frac{s^{2}}{8},\sin \frac{s^{2}}{8}\right) , \\
b &=&\left( 0,-\sin \frac{s^{2}}{8},\cos \frac{s^{2}}{8}\right) .
\end{eqnarray*}

\begin{figure}[tbp]
\centering
\subfloat[The curve $r\left(s\right) $ \label{Fig1a}] {\
\includegraphics[width=6cm,height=3cm]{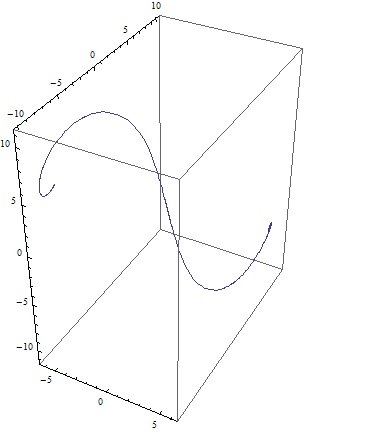}  } \hspace*{.7cm} 
\subfloat[A member of the family of surfaces having for  $a=1,b=1$ and $c=1 $ $\sigma (u)=1 .$ \label{Fig1b}] {\
\includegraphics[width=6cm,height=3cm]{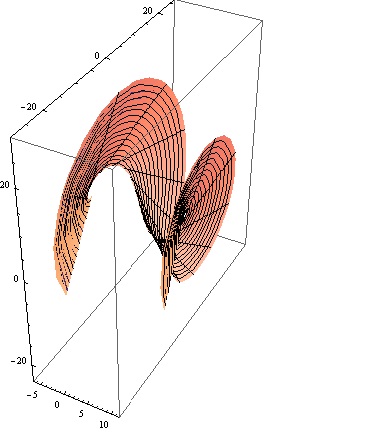}  } \newline
\par
\subfloat[A member of the family of surfaces for $a=\frac{1}{3},b=\frac{1}{5} and c=1$ $ \sigma (s)=1.$\label{Fig2a}] {\
\includegraphics[width=6cm,height=3cm]{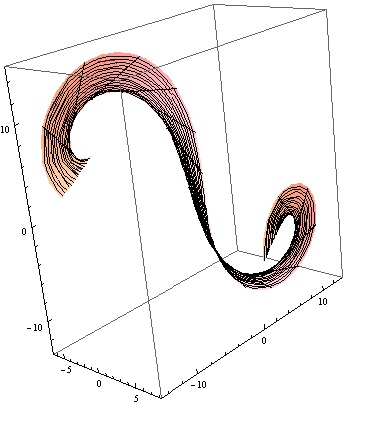}  } \hspace*{.7cm} 
\subfloat[A member of the family of surfaces for $a=\frac{1}{3},b=\frac{1}{5}$ and $c=1 $ $\sigma (s)=u $. \label{Fig2b}] {\
\includegraphics[width=6cm,height=3cm]{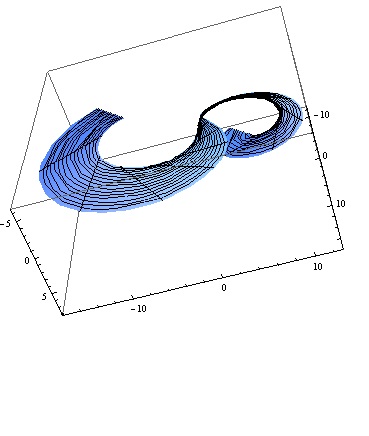}  } \newline
\par
\subfloat[The curve $r\left( s\right) $ \label{Fig3a}] {\
\includegraphics[width=6cm,height=3cm]{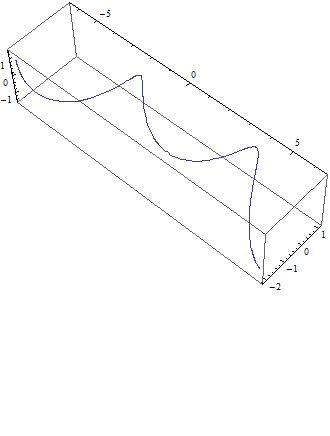}  } \hspace*{.7cm} 
\subfloat[A member of the family of surfaces for $a=1,b=1,c=1$ . \label{Fig3b}] {\
\includegraphics[width=6cm,height=3cm]{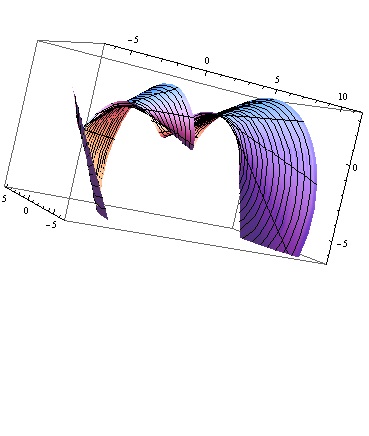}  } \newline
\par
\subfloat[A member of the family of surfaces for$a=1,b=3, c=5$.\label{Fig4a}] {\
\includegraphics[width=6cm,height=3cm]{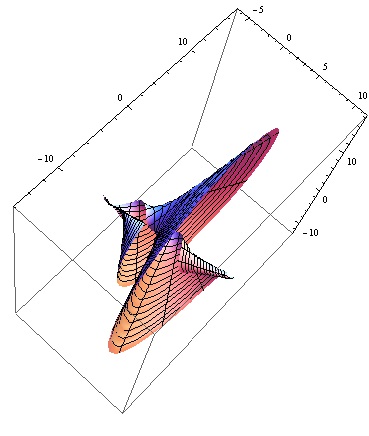}  } \hspace*{.7cm} 
\subfloat[A member of the family of surfaces for  $a=1,b=\frac{1}{5}, c=\frac{1}{10}$ . \label{Fig4b}] {\
\includegraphics[width=6cm,height=3cm]{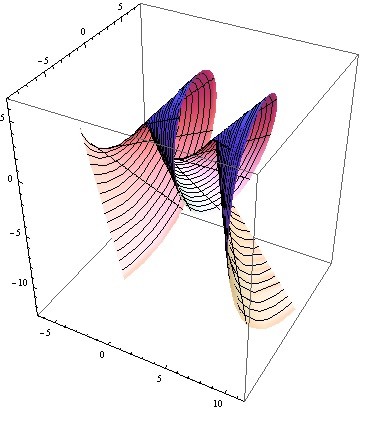}  }
\caption{}
\label{Fig2}
\end{figure}

Then, we obtain the surface pencil.

If we take $l\left( s\right) =m\left( s\right) =\sigma (s)=1,n\left(
s\right) =-1$ and $\lambda (s)=\frac{1}{2}$ 
\begin{equation*}
X\left( v\right) =v\text{ },Y\left( v\right) =\frac{3\sqrt{7}}{8}v\text{ and 
}Z\left( v\right) =\frac{1}{8}v.
\end{equation*}
Thus, a member of this family is obtained by

\begin{equation}
\varphi \left( s,v\right) =r(s)+[\alpha \left( s,v\right) t\left( s\right)
+\beta \left( s,v\right) n\left( s\right) +\gamma \left( s,v\right) b\left(
s\right) ],-2\pi \leq s\leq 2\pi ,\text{ }0\leq v\leq 5.  \label{18}
\end{equation}

Then, we plot the surface 
%TCIMACRO{\TeXButton{18}{\eqref{18}} }%
%BeginExpansion
\eqref{18}
%EndExpansion
in Fig 1b.

On the other hand, some control coefficients can be added to the function \ $%
X\left( v\right) ,Y\left( v\right) $ and $Z\left( v\right) $ such as 
\begin{equation*}
X\left( v\right) =av\text{ },Y\left( v\right) =b\frac{3\sqrt{7}}{8}v\text{
and }Z\left( v\right) =c\frac{1}{8}v,
\end{equation*}%
where $a,b$ and $c$ are real constants.

Considering $a=\frac{1}{3},b=\frac{1}{5}$ and $c=1,$ the plot of surface
pencil between the same intervals is given by Fig.1c. We obtain the shape
for taking same values of $a,b$ and $c$ and $\sigma (s)=s$ such as 
\begin{equation*}
X\left( v\right) =\frac{1}{3}v\text{ },Y\left( v\right) =\frac{1}{5}\sqrt{%
s^{2}-\frac{1}{64}}v\text{ and }Z\left( v\right) =\frac{1}{8}v,
\end{equation*}%
given by Fig. 1d.

\textbf{Example 3.9 }Let $r(s)$ be an anti-Salkowski curve given by
parametrization

\begin{equation*}
r\left( s\right) =\left( s,\frac{16}{289}[8\sin s\sinh \frac{s}{4}-15\cos
s\cosh \frac{s}{4}],-\frac{16}{289}[8\cos s\sinh \frac{s}{4}+15\sin s\cosh 
\frac{s}{4}]\right) ,
\end{equation*}%
\cite{ali} . The shape of the curve is plotted by Fig 1e. It is easy to
calculate that $\kappa =\cosh \frac{s}{4}$ and $\tau =1,$%
\begin{eqnarray*}
t &=&\left( s,\frac{16}{289}[17\sin s\cosh \frac{s}{4}+\frac{17}{4}\cos
s\sinh \frac{s}{4}],-\frac{16}{289}[17\cos s\cosh \frac{s}{4}-\frac{17}{4}%
\sin s\sinh \frac{s}{4}]\right) , \\
n &=&\left( 0,\cos s,\sin s\right) , \\
b &=&\left( 0,-\sin s,\cos s\right) .
\end{eqnarray*}

Then, the shape of the surface pencil is given by 
\begin{equation}
\varphi \left( s,v\right) =r(s)+[\alpha \left( s,v\right) t\left( s\right)
+\beta \left( s,v\right) n\left( s\right) +\gamma \left( s,v\right) b\left(
s\right) ],-2\pi \leq s\leq 2\pi ,\text{ }0\leq v\leq 5.  \label{19}
\end{equation}%
If we take $\sigma (s)=\frac{1}{\left\Vert r^{\prime \prime }(s)\right\Vert }%
=\frac{1}{\cosh \frac{s}{4}}$ and $l\left( s\right) =m\left( s\right)
,n\left( s\right) =-1$, and $\lambda (s)=\frac{\sqrt{3}}{2}.$ 
\begin{equation*}
X\left( v\right) =v,Y\left( v\right) =\sqrt{\frac{1}{\cosh \frac{s}{4}}-%
\frac{3}{4\left( \cosh \frac{s}{4}\right) ^{2}}}v\text{ and }Z\left(
v\right) =\frac{\sqrt{3}}{2}\cosh \frac{s}{4}v.
\end{equation*}

Then, we plot the surface 
%TCIMACRO{\TeXButton{19}{\eqref{19}} }%
%BeginExpansion
\eqref{19}
%EndExpansion
in Fig 1f. Considering the control coefficients $a=1,b=3$ and $c=5$, the
shape of the surface pencil between the same intervals is plotted by Figure
1g. By taking the $a=1,b=\frac{1}{5}$ and $c=\frac{1}{10},$ the plot of
surface pencil is given by Fig.1h.

\end{document}